\documentclass[a4paper,12pt]{article}
\usepackage{amsmath}
\newtheorem{theorem}{Theorem}[section]

\usepackage{amssymb}
\usepackage{amssymb,amsmath}
\usepackage{graphicx}
\usepackage{url}

\numberwithin{equation}{section}
\title{Parametric estimation for the standard and geometric telegraph process observed at discrete times}
\author{Alessandro De Gregorio\\ Department of Statistical Sciences\\ University of Padua\\ Via C.Battisti, 241, 35121 Padua - Italy  \and Stefano Maria Iacus\footnote{correspondig author, email: \texttt{stefano.iacus@unimi.it}}\\ Department of Economics, Business and Statistics\\ Via Conservatorio 7, 20122 Milan - Italy}

\begin{document}
\maketitle

\begin{abstract}
The telegraph process $X(t)$, $t>0$,  (Goldstein, 1951) and the geometric telegraph process
 $S(t) = s_0 \exp\{(\mu -\frac12\sigma^2)t + \sigma X(t)\}$ with $\mu$ a known constant and $\sigma>0$ a parameter are supposed to be  observed at $n+1$ equidistant time points
$t_i=i\Delta_n,i=0,1,\ldots, n$. For both models $\lambda$, the underlying rate of the Poisson process, is a parameter to be estimated. In the geometric case, also  $\sigma>0$ has to be estimated.
We propose different estimators of the parameters and we investigate their performance  under the high frequency  asymptotics, i.e. $\Delta_n \to 0$, $n\Delta = T<\infty$ as $n \to \infty$, with $T>0$ fixed.
The process $X(t)$ in non markovian, non stationary and not ergodic thus we use approximation arguments to derive estimators. Given the complexity of the equations involved only estimators on the first model can be studied analytically. Therefore, we run an extensive
Monte Carlo analysis to study the performance of the proposed estimators also for small sample size $n$.
\\\\
\noindent {\bf key words:} telegraph process, discretely observed
process, inference for stochastic processes.
\end{abstract}

\section{Introduction}
The random motions with finite velocity represent an alternative
to diffusion models defined by means of stochastic differential
equations. The prototype of these models is the telegraph process
(see Goldstein, 1951 and Kac, 1974) that describes the position of a particle moving
on the real line, alternatively with constant velocity $+ v$ or
$-v$. The changes of direction are governed by an homogeneous
Poisson process with rate $\lambda
>0.$ The telegraph process\footnote{In the literature, this process is alternatively called the telegraph process or the {\it telegrapher's} process.} is defined by
\begin{equation}\label{1.1}
X(t)=V(0)\int_0^t (-1)^{N(s)}{\rm d}s,\quad t>0,
\end{equation}
where the initial velocity $V(0)$ assumes the values $\pm v$ with equal probability and $\{N(t), t>0\}$ is the Poisson process on $[0,t]$. We consider a particle initially located on the real line at point $x_0$.

The explicit conditional density function of the process $X(t)$ has been obtained by
Orsingher (1990) and reads
\begin{eqnarray}
&&p(x, t;x_0,0)\label{1.2} \\
&&=\frac{e^{-\lambda t}}{2v} \left\{ \lambda
I_0\left(\frac{\lambda}{v}\sqrt{v^2\,t^2 -( x-x_0 )^2} \right)
+\frac{\partial}{\partial t} I_0\left(\frac{\lambda}{v}
\sqrt{v^2\,t^2 -( x-x_0 )^2} \right)
\right\}\chi_{\{|x-x_0|<vt\}}\notag\\
&&+\frac{e^{-\lambda t}}{2} \left\{ \delta(x-x_0-vt) +
\delta(x-x_0+vt)
\right\}\notag\\
&&= \frac{e^{-\lambda t}}{2v} \left\{ \lambda
I_0\left(\frac{\lambda}{v}\sqrt{v^2\,t^2 -( x-x_0 )^2} \right) +
\frac{v \lambda t I_1\left(\frac{\lambda}{v} \sqrt{v^2\,t^2 -(
x-x_0 )^2} \right)}{\sqrt{v^2 t^2 -( x-x_0 )^2}}
\right\}\chi_{\{|x-x_0|<vt\}}\notag\\
&&+\frac{e^{-\lambda t}}{2} \left\{ \delta(x-x_0-vt) +
\delta(x-x_0+vt)
\right\}\notag\\
&&= \frac{e^{-\lambda t}}{2v} \left\{ \lambda
I_0\left(\frac{\lambda}{v}\sqrt{u_t(x,x_0)} \right) + \frac{v
\lambda t I_1\left(\frac{\lambda}{v} \sqrt{u_t(x,x_0)}
\right)}{\sqrt{u_t(x,x_0)}}
\right\}\chi_{\{u_t(x,x_0)>0\}}+\frac{e^{-\lambda t}}{2}
\delta(u_t(x,x_0)) \notag
\end{eqnarray}
where $\{x:|x-x_0|\leq vt\}$, $u_t(x,x_0) = v^2\,t^2 -( x-x_0
)^2$ and
$$
I_\nu(x) = \sum_{k=0}^\infty \left(\frac{x}{2}\right)^{2k+\nu} \frac{1}{k!\Gamma(k+v+1)},\quad |x| < \infty, \quad |\arg x| < \pi,
$$
is the modified Bessel function with imaginary argument.
Note that the second term in equation \eqref{1.2} represents
the singular component of the distribution of \eqref{1.1}, of the position
of the particle at time $t$. Indeed, if no Poisson events occur
in the interval $[0,t]$, we have that $P\left\{X(t)=+
vt\right\}=P\left\{X(t)=-vt\right\}=\frac12 e^{-\lambda t}$.

Many authors analyzed over the years the telegraph process, see
for example Orsingher (1985, 1990), Foong and Kanno (1994), Stadje and Zacks (2004).
Di Crescenzo and Pellerey (2002) proposed the geometric telegraph
process as a model to describe the dynamics of the price of risky
assets $S(t)$. In the Black-Scholes (1973) - Merton (1973) model the process $S(t)$ is described by
means of geometric Brownian motion
\begin{eqnarray}\label{1.3}
S(t)=s_0\exp\{\alpha t+\sigma W(t)\},\quad t>0.
\end{eqnarray}
where $W(t)$ is a standard Brownian motion and
$\alpha=\mu-\frac{1}{2}\sigma^2$. Di Crescenzo and Pellerey (2002) assume that $S(t)$
evolves in time according to the following process
\begin{eqnarray}\label{1.3a}
S(t)=s_0\exp\{\alpha t+\sigma X(t)\},\quad t>0.
\end{eqnarray}
where $X(t)$ is the telegraph process. Given that $X(t)$ has bounded variation, so is $S(t)$ in equation\eqref{1.3a}. This seems a realistic way to model paths of assets in the financial markets.
Mazza and Rulliere (2004) linked the process \eqref{1.1} and the
ruin processes in the context of risk theory. Di Masi {\it et al} (1994) propose to model the volatility of financial markets in terms of the telegraph process. Ratanov (2004, 2005) propose to model financial markets using a telegraph process with two intensities $\lambda_\pm$ and two velocities $c_\pm$. While such markets may admit an arbitrage opportunity, linking opportunity velocities and interest rates, the author proves that the market becomes arbitrage-free and complete. An analogous of the Black and Scholes equation is established as well.

The aim of this paper is the estimation of the parameter $\lambda$
when $\{X(t), 0 \leq t \leq T\}$ is observed at
equidistant times $0=t_0<...<t_n$ (and also $\sigma$ for discrete observations from the process \eqref{1.3a}). We assume that
$t_i=i \Delta_n$, $i=0, \ldots, n$, hence $n\Delta_n=T$. The
asymptotic framework is the following: $n\Delta_n=T$ fixed and
$\Delta_n \rightarrow 0$ as $n\to \infty$. Sometimes we will use $\Delta$ instead of $\Delta_n$ to simplify the writing.
When the telegraph process $X(t)$ is observed continuously then $N(T)/T$ is the optimal estimator of the parameter $\lambda$ as this statistical problem will be equivalent to the one of the observation of the whole  Poisson process\footnote{For more details on parametric estimation for Poisson process see Kutoyants (1998).} on $[0,T]$. This situation corresponds indeed to the limiting  experiment in our asymptotic framework.
This asymptotic framework is usually known as the ``high frequency'' scheme in the literature on estimation from discrete time observations of processes solution to stochastic differential equations of the form
\begin{eqnarray}
d Y(t)=b(Y(t),\theta)dt+\sigma(Y(t),\theta)dW(t),
\end{eqnarray}
This field has been an active research area during the last twenty years. The
reader can consult Sorensen (2004) for a review on  estimation techniques
recently appeared in the literature since the seminal papers of Le Breton (1976) and Florens-Zmirou (1989): e.g. estimating functions, analytical and
numerical approximations of the likelihood function, MCMC methods,
indirect inference, etc. Unfortunately such methods are not directly applicable in our case because the telegraph process is not ergodic or stationary nor Markovian.
The main idea in the paper is to consider the observed increments of the process $X(i\Delta) - X((i-1)\Delta)$ as $n$ copies of the telegraph process up to time $\Delta$ and treat them as if they were independent (which is  untrue). From this idea we build an approximated likelihood and score function from which we derive estimators. We further propose least squares estimators. Equations emerging in connection with the telegraph process are always complicated to treat and closed form results are quite rare in the literature. This also happens in our work, thus in some cases we rely on numerical simulations to study the properties of the estimators.

It is worth mention that, up to our knowledge, the only references about estimation problems for the telegrapher's  processes are  Yao (1985) and Iacus (2001). The first author considers a the problem of state estimation of the telegrapher's process under white noise perturbation and studies performance of nonlinear filters. The second paper is about the estimation of the parameter $\theta$ of the non-constant rate $\lambda_\theta(t)$ from continuous observations of the process.

The paper is organized as follows. In Section \ref{sec2} we will
introduce the approximating likelihood function
 \begin{eqnarray}
L_n(\lambda ) =& \prod_{i=1}^n p(X_i, \Delta_n;X_{i-1},t_{i-1})
\end{eqnarray}
where $p(X_i, \Delta_n;X_{i-1},t_{i-1})$ is defined by
\eqref{1.2}. We will then study the asymptotic properties of
estimator 
$$\bar\lambda_n = \arg\max_{\lambda>0} L_n(\lambda )$$ 
(which is interpreted as an approximated maximum likelihood estimator) and the properties of
the estimator obtained on the approximated score function whose properties are easier to study, i.e. the quantity $\hat\lambda_n$ satisfying
$$ \left.\frac{\partial}{\partial \lambda} \log L_n(\lambda )\right|_{\lambda=\hat\lambda_n} = 0$$ 
We also propose a least squares estimator for $\lambda$ based on the second moment of the process $X(t)$. Inference problems and estimators for the parameters $\lambda$ and $\sigma$ of the geometric telegrapher's process are considered in Section \ref{sec:geo}.
Finally, section \ref{sec3} contains a Monte Carlo analysis  to study
 empirically the behavior of the estimators  in a finite sample context (i.e. non asymptotically).
 
\section{The scheme of observation and the asymptotics}\label{sec2}
We assume that the telegraph process $\{X(t), 0 \leq t \leq T\}$,
with $X(0)=x_0=0,$ is observed only at discrete times $0 <
t_1 < \cdots <t_n=T$, with $t_i = i\Delta_n$, $i=0,
\ldots, n$ hence $n\Delta_n = T$. We use the following notation to
simplify the formulas: $X(t_i) = X(i\Delta_n) = X_i$. The
asymptotic is considered as $n$ tends to infinity under the
conditions $\Delta_n\to 0$ and $n\Delta_n = T$. The interest is in
the estimation of the parameter $\lambda$ whilst $v$ is assumed to
be known.

As mentioned in the Introduction, if one can observe the whole trajectory, $\lambda$ can be
estimated as $N(T)/T$ where $N(T)$ is the number of times the
process switches its velocity during the interval $[0,T]$ which
is, of course, the number of Poisson events counted in $[0,T]$.
This is certainly the best estimator of $\lambda$ and it is indeed our target.

The estimation of $v$ is always an uninteresting problem as, if
there are no switchings in $[(i-1)\Delta_n, i\Delta_n]$ then $X_i
- X_{i-1} = v \Delta_n$, hence if $\Delta_n$ is sufficiently
small, there is high probability of observing $N(t_{i+1}) - N(t_i)
= 0$ then $v$ can be estimated (actually calculated) without
error.

The process $X(t)$ itself is not Markovian. On the contrary, the
two dimensional process $(X(t), V(t))$ has the Markov property but
 a scheme of observation in which one is able
to observe both the position and the velocity of the process at
discrete time instants is not admissible, so we can rely only on the observation of
the $X(t)$ component. Hence we cannot write an explicit likelihood
of the process in the form of a product of transition densities as, for example, for the case of diffusion processes. Another unfortunate fact about the telegraph process is that it is not even stationary  (at the
second order) as
\begin{equation}
E X(t) = 0
\label{eq:m1}
\end{equation}
and
\begin{equation}
E X^2(t)
=\frac{v^2}{\lambda}\left(t-\frac{1-e^{-2\lambda
t}}{2\lambda}\right)
\label{eq:m2}
\end{equation}
(see Orsingher (1990)) nor it posseses an ergodic property, so
we cannot use the same approach as proposed, e.g.,  in Sorensen (2000). We then need an approximation argument as follows.
\subsection{An approximation of the likelihood function}
By taking into account the distribution \eqref{1.2}, we
approximate the likelihood of the process with the following
function
\begin{eqnarray}
L_n(\lambda)&=&L_n(\lambda | X_0, X_1, \ldots, X_n) =
\prod_{i=1}^n
p(X_i, \Delta_n;X_{i-1},t_{i-1})\label{2.2}\\
&=&\prod_{i=1}^n\Bigg\{ \frac{e^{-\lambda\Delta_n}}{2v} \left\{
\lambda I_0\left(\frac{\lambda}{v}\sqrt{u_{n,i}} \right) + \frac{v
\lambda \Delta_n I_1\left(\frac{\lambda}{v} \sqrt{u_{n,i}}
\right)}{\sqrt{u_{n,i}}}
\right\}\chi_{\{u_{n,i}>0\}}\notag\\
&&+\frac{e^{-\lambda \Delta_n}}{2} \delta(u_{n,i}=0)\Bigg\}\notag
\end{eqnarray}
where $u_{n,i} =  u_n(X_i,X_{i-1}) = v^2 \Delta_n^2 -  (
X_i-X_{i-1} )^2$.

The density $p(X_i, \Delta_n;X_{i-1},t_{i-1})$ appearing in
\eqref{2.2} is  the probability law
of a telegraph process initially located in $X_{i-1}$, that
reaches the position $X_{i}$ at time $t_i$. To build our approximation of the likelihood,  we consider the observed increments $X_i - X_{i-1}$ as $n$ copies of the process $X(\Delta)$ (i.e. the process $X(t)$ up to time $\Delta$)  and we treat them as if they were independent. This is of course untrue, but results below show that our is not a bad idea.
 It is clear that \eqref{2.2} is
equivalent to
\begin{eqnarray}
L_n(\lambda)&=&\left(\frac{e^{-\lambda
\Delta_n}}{2}\right)^{n-n^+}\prod_{i=1}^{n^+}
\frac{e^{-\lambda\Delta_n}}{2v} \left\{ \lambda
I_0\left(\frac{\lambda}{v}\sqrt{u_{n,i}} \right)+ \frac{v \lambda
\Delta_n I_1\left(\frac{\lambda}{v} \sqrt{u_{n,i}}
\right)}{\sqrt{u_{n,i}}} \right\}\label{2.3}\\
&=&\frac{e^{-\lambda n
\Delta_n}}{2^n}\frac{1}{v^{n^+}}\prod_{i=1}^{n^+} \left\{ \lambda
I_0\left(\frac{\lambda}{v}\sqrt{u_{n,i}} \right)+ \frac{v \lambda
\Delta_n I_1\left(\frac{\lambda}{v} \sqrt{u_{n,i}}
\right)}{\sqrt{u_{n,i}}} \right\}\notag
\end{eqnarray}
where $n^+$ is equal to the number of approximating telegraph
processes with at least one change of direction (see Figure 1).

\begin{figure}[t]
\includegraphics{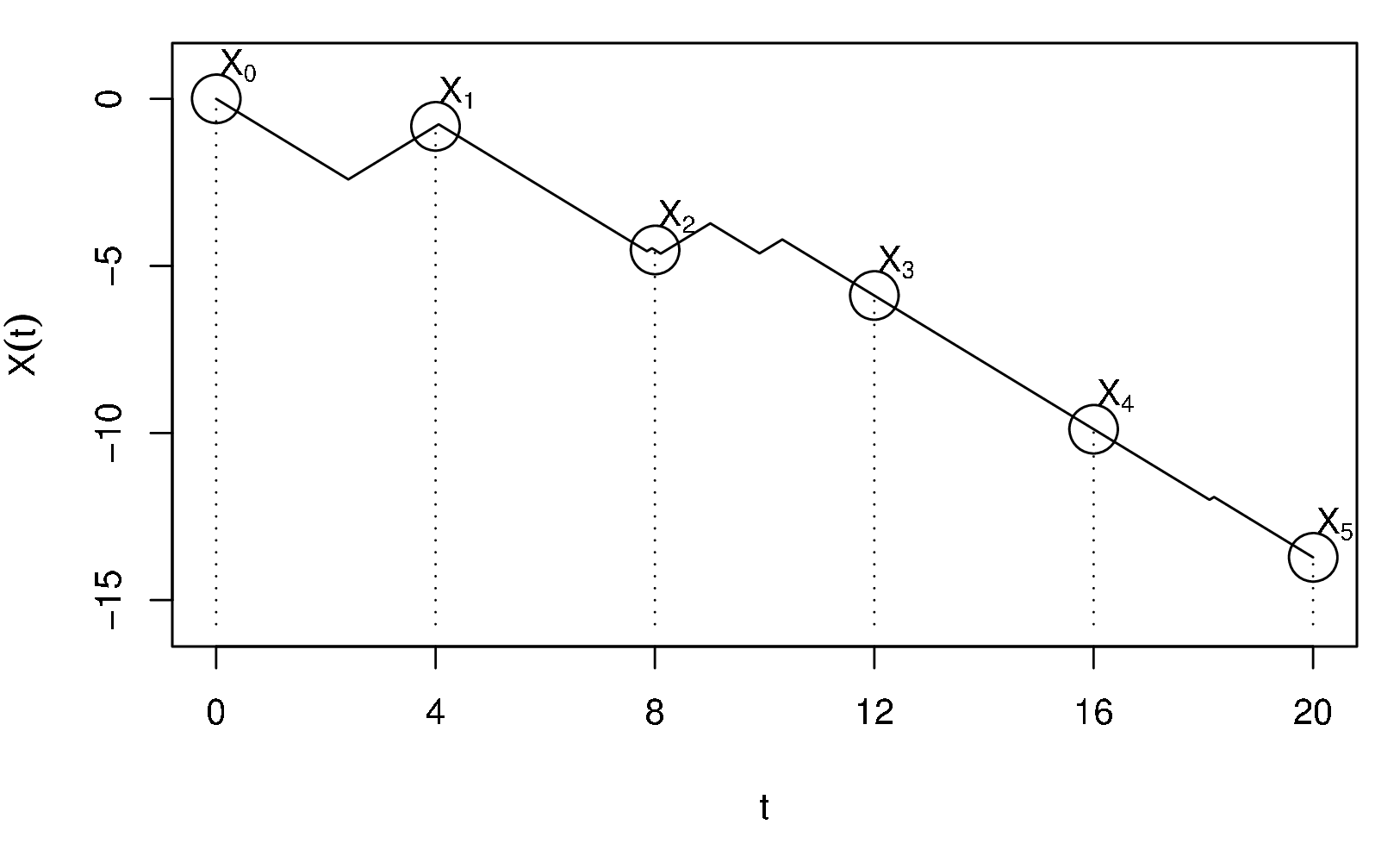}
\caption{Discrete time sampling of the telegraph process. Between
$X_3$ and $X_4$ no Poisson event occurred with probability
$e^{\lambda\Delta_n}$, hence the singular part of the approximated
likelihood emerges in the $L_n(\lambda)$. For this example trajectory $n=5$ and $n^+ =4$.}
\label{fig1}
\end{figure}

 In
the expression \eqref{2.3}, the factor $\left(\frac{e^{-\lambda
\Delta_n}}{2}\right)^{n-n^+}$ concerns the singular part of the
densities $p(X_i, \Delta_n;X_{i-1},t_{i-1})$, while the product
represents the absolutely continuous components of the
distributions of the telegraph processes. Note that for
increasing values of $\lambda$, the absolutely continuous
component of \eqref{2.3} has a bigger weight than the discrete
component; viceversa for small values of $\lambda$ this has consequences in the performance of the estimators as shown in Section \ref{sec3}. Figure \ref{fig1} shows how the two components of the function $L_n(\lambda)$ emerge for this scheme of observation.
\subsection{Estimators on the approximated likelihood}
One estimator that can be derived from $L_n(\lambda)$ is the following quantity
\begin{eqnarray}\label{2.4}
\bar\lambda_n = \arg\max_{\lambda>0} L_n(\lambda)
\end{eqnarray}
Of course, $\bar\lambda_n$  is not a true maximum likelihood estimator as $L_n(\lambda)$ is not itself a true likelihood, nevertheless such estimators can be effective like in Kessler (2000) but unfortunately, our model is not ergodic or stationary and the likelihood  looks  quite difficult to handle. Nevertheless we are able to prove uniqueness of the estimator.
On the contrary, it is not clear how to prove asymptotic properties of the estimator directly on  \eqref{2.4} as the term $\log(2^n v^{n^+})$ diverges in
\begin{equation}
\label{logL}
\log L_n(\lambda)=-\lambda n \Delta_n
-\log(2^n v^{n^+})+
\sum_{i=1}^{n^+}\log \left\{ \lambda
I_0\left(\frac{\lambda}{v}\sqrt{u_{n,i}} \right)+ \frac{v \lambda
\Delta_n I_1\left(\frac{\lambda}{v} \sqrt{u_{n,i}}
\right)}{\sqrt{u_{n,i}}} \right\}
\end{equation}
On the contrary the following estimating function 
\begin{equation}
F(\lambda; X_1, \ldots, X_n) = \frac{\partial}{\partial \lambda} \log L_n(\lambda)
\label{est.fun}
\end{equation}
is easier to study. The estimator in this case is the value $\hat\lambda_n$  solution of $F(\lambda)=0$, i.e.
\begin{equation}
\hat\lambda_n : F(\lambda=\hat\lambda_n; X_1, \ldots, X_n) =0
\label{eq:lamhat}
\end{equation}
and once again, we borrow the approach from the literature on estimating functions for discretely observed diffusion processes.
The function \eqref{est.fun} is indeed more tractable and the asymptotic properties of the estimator $\hat\lambda_n$ can be studied as in the next theorem. We just note that,
given the uniqueness of $\bar\lambda_n$ (proved in Theorem \ref{th:unique} below),  $\hat\lambda_n$ coincides with $\bar\lambda_n$ hence they share the same properties and they will be treated as one in Table \ref{tab1}.
\begin{theorem}
\label{th2}
Under the condition $n\Delta_n=T$ as $\Delta_n\rightarrow 0$ we
have that
\begin{eqnarray}\label{2.5}
\hat\lambda_n\rightarrow \hat\lambda_\infty=\frac{N(T)}{T}\,.
\end{eqnarray}
\end{theorem}
\textbf{Proof.} In order to prove \eqref{2.5} we recall some properties of Bessel's functions (see e.g. \S 5.7, Lebedev, 1972)
$$
\frac{{\rm d}}{{\rm d} x} I_n(x) = \frac12\biggl( I_{n-1}(x) +  I_{n+1}(x)\biggr)
$$
and
$$
\lim_{u\to 0} \frac{I_1(k \cdot u)}{u} = \frac{k}{2},\quad
\lim_{u\to 0} I_0(k\cdot u) =1,\quad
\lim_{u\to 0} I_1(k\cdot u) =0,\quad
\lim_{u\to 0} I_2(k\cdot u) =0,
$$
Direct differentiation of \eqref{logL} gives
\begin{equation}
-n \Delta_n + \sum_{i=1}^{n^+} 
\frac{I_0\left(\frac{\lambda}{v}\sqrt{u_{n,i}}\right)\left(1+ \frac{\lambda\Delta_n}{2}\right) + 
\left( \frac{\lambda\sqrt{u_{n,i}}}{v} +  \frac{v \Delta_n}{\sqrt{u_{n,i}}}\right)I_1\left(\frac{\lambda}{v}\sqrt{u_{n,i}}\right) + 
\frac{\lambda\Delta_n}{2} I_2\left(\frac{\lambda}{v}\sqrt{u_{n,i}}\right)}{\lambda I_0\left(\frac{\lambda}{v}\sqrt{u_{n,i}}\right)
+ \frac{\lambda v\Delta_n}{\sqrt{u_{n,i}}} I_1\left(\frac{\lambda}{v}\sqrt{u_{n,i}}\right) }
\label{dlogL}
\end{equation}
In the limit as $n\to\infty$ we have, $\Delta_n\to 0$, $n\Delta_n\to T$,  $u=u_{n,i} = v^2 \Delta_n^2 -  ( X_i-X_{i-1} )^2\rightarrow
0$ and $n^+\rightarrow N(T)$. Therefore \eqref{dlogL} converges to
$$
-T + \frac{N(T)}{\lambda}\,.
$$
\hfill$\square$
\\
In light of Theorem \ref{th2} we can say that in the high frequency
observation scheme, the estimator $\hat{\lambda}_n$ tends to the best
estimator of the intensity $\lambda$ of an homogeneous Poisson
process. 
\begin{theorem}
The estimator  $\bar\lambda_n$ in  \eqref{2.4} is unique.
\label{th:unique}
\end{theorem}
\textbf{Proof.} 
To prove the result, we show that the second partial derivative of $\log L_n(\lambda)$ is negative.
In view of the following property of Bessels functions
\begin{equation}
I_{\nu-1}(z) - I_{\nu+1}(z) = \frac{2 \nu}{z} I_\nu(z)
\label{eq:bess}
\end{equation} 
we have that
$$
-\frac{\lambda\Delta_n}{2}\left\{I_0\left(\frac{\lambda}{v}\sqrt{u_{n,i}}\right) -  I_2\left(\frac{\lambda}{v}\sqrt{u_{n,i}}\right)\right\} = - \frac{\Delta_n v}{\sqrt{u_{n,i}}} I_1\left(\frac{\lambda}{v}\sqrt{u_{n,i}}\right) 
$$
therefore equation \eqref{dlogL} can be rewritten as
\begin{equation}
\frac{\partial}{\partial\lambda} \log L_n(\lambda) =  -n \Delta_n + \sum_{i=1}^{n^+} 
\frac{{\sqrt{u_{n,i}}}\,v\,\left( 1 + \Delta_n \,\lambda  \right) \,I_0\left(\frac{\lambda}{v}\sqrt{u_{n,i}}\right) + 
    u_{n,i}\,\lambda \,I_1\left(\frac{\lambda}{v}\sqrt{u_{n,i}}\right)}{v\,\lambda \,
    \left( {\sqrt{u}}\,I_0\left(\frac{\lambda}{v}\sqrt{u_{n,i}}\right) + 
      v\,\Delta_n \,I_1\left(\frac{\lambda}{v}\sqrt{u_{n,i}}\right) \right) }
      \label{eq:diff2}
\end{equation}
Differentiating \eqref{eq:diff2} again with respect to $\lambda$ and posing $x =  \frac{\lambda}{v}\sqrt{u_{n,i}}$ to simplify the equations\footnote{We remind that this $x = \sqrt{u_{n,i}}\lambda / v$ is strictly positive because it refers to the $n^+$ terms for which $u_{n,i}>0$.}, we obtain
$$
\frac{\partial^2}{\partial\lambda^2} \log L_n(\lambda) =
\sum_{i=1}^{n^+} \frac{g_i}{2\,v^2\,{\lambda }^2\,
    {\left( {\sqrt{u_{n,i}}}\,I_0(x)  + 
        v\,\Delta_n \,I_1(x)  \right) }^2}
$$
where the generic term $g_i$ is the sum of the following terms
\begin{equation}
 -2\,{\sqrt{u_{n,i}}}\,v^3\,\Delta_n \,I_1(x)   I_0(x)  <0
\label{eq3}
\end{equation}
\begin{equation}
u_{n,i}\,\left\{ u_{n,i}\,{\lambda }^2 - v^2\,\left( 2 + \Delta_n \,\lambda  + \Delta_n^2\,{\lambda }^2 \right)  \right\} \, I_0^2(x) 
\label{eq1}
\end{equation}
\begin{equation}
  2\,u_{n,i}\,\lambda \,\left\{ -\lambda u_{n,i}  + v^2\,\Delta_n \,\left( 1 + \Delta_n \,\lambda  \right) 
       \right\} \, I_1^2(x) 
\label{eq2}
\end{equation}
\begin{equation}
      u_{n,i}\,\lambda \,\left( u_{n,i}\,\lambda  - v^2\,\Delta_n \,\left( 1 + \Delta_n \,\lambda  \right)  \right) \,
        I_2(x)   I_0(x)   
\label{eq4}
\end{equation}
Recalling that $u_{n,i} = v^2 \Delta_n^2 - (X_i - X_{i-1})^2 > 0$, equation \eqref{eq1} can be rewritten as
\begin{equation}
 - u_{n,i} \,\left\{ v^2\,\Delta_n \,\lambda + (X_i - X_{i-1})^2\,{\lambda }^2 
       \right\} I_0^2(x)  - 2 u_{n,i}
    v^2\, I_0^2(x)
\label{eq1.1}
\end{equation}
 summation of   \eqref{eq2} and \eqref{eq4} gives
\begin{equation}
-u_{n,i} \, \,\left\{ v^2\,\Delta_n\lambda  + (X_i - X_{i-1})^2\,\lambda^2  \right\} \,
  \left( 2\,{I_1(x)}^2 - 
    I_0(x)    I_2(x) \right)
\label{eq4.1}
\end{equation}
Putting together
$$
 - u_{n,i} \,\left\{ v^2\,\Delta_n \,\lambda + (X_i - X_{i-1})^2\,{\lambda }^2 
       \right\} I_0^2(x)
$$
and equation \eqref{eq4.1} it remains to study the sign of
$$
I_0^2(x) + 2\,{I_1(x)}^2 -    I_0(x)    I_2(x) =
I_0(x)(I_0(x) - I_2(x)) + 2\,{I_1(x)}^2
$$
which is positive due to the fact that $I_0(x) > I_2(x)$ for positive $x$ (from property \eqref{eq:bess}).
\hfill$\square$\\
\subsection{A least squares estimator}
We already mentioned that the telegrapher's process is not stationary at the second order, but we can still think to use \eqref{eq:m2} to obtain a least squares estimator on the mean of the squared increments in the following way.
If, as before, we consider the observations $X_i-X_{i-1}$ as $n$ copies of the process $X(\Delta)$, then  looking at \eqref{eq:m2} for $t=\Delta$ we have
$$
E X^2(\Delta)
=\frac{v^2}{\lambda}\left(\Delta-\frac{1-e^{-2\lambda
\Delta}}{2\lambda}\right)
$$
Consider now the sample second moment (or the mean of the square) of the observed increments $X_i - X_{i-1}$
$$
m_2 = \frac{1}{n} \sum_{i=1}^n (X_i - X_{i-1})^2
$$
then the following estimator can be considered
\begin{equation}
\tilde\lambda_n = \arg\min_{\lambda>0} \left\{m_2 - \frac{v^2}{\lambda}\left(\Delta_n-\frac{1-e^{-2\lambda
\Delta_n}}{2\lambda}\right)\right\}^2
\label{eq:lamtil}
\end{equation}
This estimator will be compared numerically to the estimator $\hat\lambda_n$ in Section \ref{sec3}.

\section{Parametric estimation for the geometric telegraph process}\label{sec:geo}
Consider the process $Y$ of the observed log-returns
$$
Y_i = \log \frac{S_i}{S_{i-1}} = \alpha \Delta_n + \sigma (X_i - X_{i-1})
$$
where $S_i = S(t_i)$ are discrete observations from the geometric telegraph process \eqref{1.3a}.
We assume $\mu$ to be known, which is usually the case in finance where $\mu$ is related to the expected return of non risky assets like bonds, etc. The parameters $\sigma$ and $\lambda$ are to be estimated. As in the previous sections, we can assume $v$ to be known as well, if not we will show in the next paragraph a simple way to obtain it.
We assume $Y_i$ to be $n$ copies of the process 
$$Y(\Delta) = \alpha \Delta + \sigma X(\Delta) $$
with $X(\Delta) = X_i - X_{i-1}$ and $X(0) = x_0 = 0$. Therefore, by \eqref{eq:m1}, we have
$${\rm E} Y(\Delta) = \alpha \Delta$$
and by \eqref{eq:m2} we obtain
\begin{equation}
{\rm Var} Y(\Delta) = \sigma^2 {\rm Var} X(\Delta) = \sigma^2 \frac{v^2}{\lambda}\left(
\Delta - \frac{1-e^{-2\lambda\Delta}}{2\lambda}
\right)
\label{eq:v2}
\end{equation}
A good estimator of the volatility $\sigma$ can be derived from  the sample mean of the log returns.
Indeed,
\begin{equation}
\bar Y_n = \frac{1}{n} \sum_{i=1}^n Y_i = \alpha \Delta + \frac{\sigma}{n} X_n
\label{eq:y1}
\end{equation}
and
\begin{equation}
{\rm E} \bar Y_n = \alpha \Delta + \frac{\sigma}{n} {\rm E} X_n = \alpha \Delta = \left(\mu - \frac12 \sigma^2\right) \Delta
\label{eq:y2}
\end{equation}
again by \eqref{eq:m1} and for the properties of the log-returns. 
From \eqref{eq:y2} we have that
$$
\sigma^2 = 2\left(\mu - \frac{{\rm E} \bar Y_n}{\Delta}\right)
$$
from which the following unbiased estimator of $\sigma^2$ can be derived
$$
\hat \sigma_n^2 = 2\left(\mu - \frac{\bar Y_n}{\Delta}\right)\,.
$$
Therefore, a reasonable estimator of $\sigma$ is
\begin{equation}
\hat \sigma_n = \sqrt{2\left(\mu - \frac{\bar Y_n}{\Delta}\right)}
\label{eq:sigma}
\end{equation}
which not always exists because there is
no guarantee that $\mu > \bar Y_n/\Delta$. Moreover, it should be noticed that
in practice, given $\mu$, $\sigma$ and $\Delta$ the estimator essentially depends on the last value of the telegraph process $X_n$. In fact, we can write \eqref{eq:sigma} in terms of the telegraph process
$$
\sqrt{2\left(\mu - \frac{\bar Y_n}{\Delta}\right)} = 
 \sqrt{2\left(\mu - \alpha - \frac{\sigma}{T} X_n\right)}
$$
therefore the estimate of $\sigma$ does not depend on $n$. This is why Table \ref{tab4} reports the same value of the estimates for different sample sizes.
We then use $\hat\sigma_n$ to estimate $\lambda$ making use of \eqref{eq:v2}.
Let 
$$
\bar s^2_Y = \frac{1}{n} \sum_{i=1}^n (Y_i-\bar Y_n)^2
$$
then the proposed estimator of $\lambda$ is 
\begin{equation}
\dot\lambda_n = \arg\min_{\lambda>0} \left( \bar s^2_Y -  \hat\sigma^2_n \frac{v^2}{\lambda}\left(
\Delta - \frac{1-e^{-2\lambda\Delta}}{2\lambda}
\right)\right)^2
\label{eq:lamdot}
\end{equation}
\subsection{Filtering of the geometric telegraph process}
If the velocity $v$ is not known one can proceed as follows: set
$$
Z_i = \frac{Y_i - {\rm E}\bar Y_n}{\sigma} = X_i - X_{i-1} = X(\Delta)
$$
hence an estimator of the increments of the telegraph process is
$$
\hat Z_i = Y_i - \frac{\bar Y_n}{\hat\sigma_n} = \hat X(\Delta),\quad i=1, \ldots, n
$$
then
$$
\hat Z_1 = \hat X_1,\qquad \hat Z_2 + \hat Z_1 = \hat X_2, \qquad \hat Z_3 + \hat Z_2 + \hat Z_1 = \hat X_3, \ldots
$$  
where $\hat X_i$ are the estimated states of the underlying telegrapher's process. From these estimates, one can proceed as in previous sections and estimate both $\lambda$ and $v$.

\section{Monte Carlo analysis}\label{sec3}
  To assess the properties of the estimator  \eqref{2.4} and \eqref{eq:lamhat} for fixed $n<\infty$ we run extensive Monte Carlo
  analysis. Given that numerically $\bar\lambda_n$ and $\hat\lambda_n$ coincides, we treat them as one in the tables.
 We simulate 10000 trajectories of the telegrapher's process on the interval $[0,T]$, $T=500$, for different values of $\lambda$ and with $v=1$ fixed. Each trajectory has then been resampled on a regular grid of $n = 50$, 100, 500 and 1000 points and the corresponding observations have been used to estimate the unknown parameter. The results have been collected in Table \ref{tab1}. It emerges that, as expected, the bias tends asymptotically to zero as well as the mean square error. Furthermore, bias and variance are strictly correlated to the value of the unknown parameter $\lambda$. This is expected as well because, for fixed $n$, as the more $\lambda$ increases the more Poisson events remain hidden to the observer. For the same experiment (and on the same sample trajectories),  Table \ref{tab2} reports the performance on the least squares estimator $\tilde\lambda_n$ from equation \eqref{eq:lamtil}.

Tables \ref{tab3} and \ref{tab4} reports estimates results on the telegrapher's process respectively for the estimation of $\lambda$ and $\sigma$. The paths of the geometric telegrapher's process have been generated from the ones of the telegrapher process of Tables \ref{tab1} and \ref{tab2}.
As it can be seen, the estimator $\dot\lambda_n$ in Table \ref{tab3} strongly depends on the quality of the estimate $\hat\sigma_n$ (reported in Table \ref{tab4}). For low values of $\lambda$, in some cases the condition for the existence of $\hat\sigma_n$, i.e.   $\mu > \bar Y_n/\Delta$, has not been fulfilled hence we report the percentage of valid paths over the 10000 simulated.
For $\lambda\geq 0.75$, it seems that $\dot\lambda_n$ performs quite similarly to $\tilde\lambda_n$ in terms of bias. This seems consistent with the definition of the estimators and the performance of the estimator $\hat\sigma_n$.
Tables also report the column $\sqrt{{\rm MSE}}$. Values under this column are calculated, e.g., as follows
$$
\sqrt{{\rm MSE(\hat\lambda_n)}} =
\sqrt{\frac{1}{N} \sum_{i=1}^N (\hat\lambda_n - \lambda)^2}
$$
where $N$ is the number of Monte Carlo simulations ($N=10000$ in our case) and $n$ is the fixed sample size.
% latex table generated in R 2.3.1 by xtable 1.3-2 package
% Thu Jul  6 00:10:52 2006
\begin{table}[ht]
\begin{center}
\begin{tabular}{r|rrrr|r}
\hline
 $\lambda$ & Bias & $\sqrt{{\rm MSE}(\lambda)}$ & $\min\hat\lambda$ & $\max\hat\lambda$ & $n$ \\
\hline
 0.10 & $-$0.002 & 0.018 & 0.04 & 0.20 &  50 \\
 & $-$0.001 & 0.016 & 0.05 & 0.16 & 100 \\
 & $-$0.000 & 0.014 & 0.06 & 0.16 & 500 \\
 & $-$0.000 & 0.014 & 0.05 & 0.15 & 1000 \\
 0.25 & $-$0.011 & 0.041 & 0.13 & 0.47 &  50 \\
 & $-$0.003 & 0.031 & 0.16 & 0.41 & 100 \\
 & $-$0.000 & 0.023 & 0.16 & 0.34 & 500 \\
 & $-$0.000 & 0.023 & 0.16 & 0.35 & 1000 \\
 0.50 & $-$0.062 & 0.092 & 0.26 & 0.78 &  50 \\
 & $-$0.011 & 0.059 & 0.32 & 0.85 & 100 \\
 & $-$0.001 & 0.035 & 0.37 & 0.65 & 500 \\
 & $-$0.000 & 0.033 & 0.37 & 0.63 & 1000 \\
 0.75 & $-$0.151 & 0.175 & 0.36 & 1.01 &  50 \\
 & $-$0.031 & 0.091 & 0.47 & 1.18 & 100 \\
 & $-$0.001 & 0.048 & 0.58 & 0.96 & 500 \\
 & $-$0.000 & 0.043 & 0.60 & 0.92 & 1000 \\
 1.00 & $-$0.264 & 0.283 & 0.45 & 1.23 &  50 \\
& $-$0.064 & 0.128 & 0.62 & 1.53 & 100 \\
 & $-$0.001 & 0.058 & 0.79 & 1.26 & 500 \\
 & $-$0.001 & 0.051 & 0.81 & 1.22 & 1000 \\
 1.50 & $-$0.546 & 0.558 & 0.58 & 1.48 &  50 \\
 & $-$0.162 & 0.227 & 0.90 & 2.16 & 100 \\
 & $-$0.003 & 0.080 & 1.24 & 1.86 & 500 \\
 & $-$0.001 & 0.066 & 1.27 & 1.77 & 1000 \\
 2.00 & $-$0.874 & 0.882 & 0.75 & 1.65 &  50 \\
 & $-$0.298 & 0.357 & 1.11 & 2.67 & 100 \\
 & $-$0.006 & 0.106 & 1.63 & 2.47 & 500 \\
 & $-$0.000 & 0.083 & 1.68 & 2.33 & 1000 \\
\hline\hline
\end{tabular}
\end{center}
\caption{Empirical performance of the estimator $\hat\lambda_n$ of \eqref{2.4} for
different values of the parameter $\lambda$ and different sample
size. The time horizon $T$ is fixed to 500. Results over 10000
Monte Carlo paths of the telegraph process.  See text for more details.} 
\label{tab1}
\end{table}

\begin{table}[ht]
\begin{center}
\begin{tabular}{r|rrrr|r}
\hline
 $\lambda$ & Bias & $\sqrt{{\rm MSE}(\lambda)}$ & $\min\hat\lambda$ & $\max\hat\lambda$ & $n$ \\
\hline
 0.10 & 0.002 & 0.022 & 0.04 & 0.24 &  50 \\
 & 0.001 & 0.018 & 0.05 & 0.18 & 100 \\
 & 0.000 & 0.016 & 0.05 & 0.16 & 500 \\
 & $-$0.000 & 0.016 & 0.05 & 0.16 & 1000 \\
 0.25 & 0.007 & 0.051 & 0.13 & 0.55 &  50 \\
 & 0.003 & 0.037 & 0.14 & 0.41 & 100 \\
 & 0.000 & 0.026 & 0.16 & 0.36 & 500 \\
 & 0.000 & 0.025 & 0.17 & 0.35 & 1000 \\
 0.50 & 0.018 & 0.106 & 0.25 & 1.13 &  50 \\
 & 0.007 & 0.070 & 0.32 & 0.93 & 100 \\
 & 0.001 & 0.040 & 0.35 & 0.66 & 500 \\
 & 0.000 & 0.037 & 0.35 & 0.65 & 1000 \\
 0.75 & 0.028 & 0.161 & 0.40 & 1.72 &  50 \\
 & 0.012 & 0.105 & 0.45 & 1.32 & 100 \\
 & 0.001 & 0.054 & 0.56 & 1.01 & 500 \\
 & 0.000 & 0.048 & 0.58 & 0.96 & 1000 \\
 1.00 & 0.040 & 0.219 & 0.53 & 2.42 &  50 \\
 & 0.017 & 0.141 & 0.62 & 1.80 & 100 \\
 & 0.002 & 0.066 & 0.73 & 1.27 & 500 \\
 & 0.001 & 0.057 & 0.77 & 1.25 & 1000 \\
 1.50 & 0.059 & 0.329 & 0.72 & 3.00 &  50 \\
 & 0.028 & 0.218 & 0.92 & 2.78 & 100 \\
 & 0.003 & 0.093 & 1.19 & 1.89 & 500 \\
 & 0.001 & 0.075 & 1.25 & 1.79 & 1000 \\
 2.00 & 0.080 & 0.412 & 1.05 & 3.00 &  50 \\
 & 0.035 & 0.290 & 1.22 & 3.00 & 100 \\
 & 0.006 & 0.120 & 1.59 & 2.46 & 500 \\
 & 0.003 & 0.095 & 1.66 & 2.39 & 1000 \\
\hline\hline
\end{tabular}
\end{center}
\caption{Empirical performance of the estimator $\tilde\lambda_n$ of \eqref{eq:lamtil} for
different values of the parameter $\lambda$ and different sample
size. The time horizon $T$ is fixed to 500. Results over 10000
Monte Carlo paths of the telegraph process. See text for more details.} 
\label{tab2}
\end{table}

\begin{table}[ht]
\begin{center}
\begin{tabular}{r|rrrr|rr}
\hline
  $\lambda$ & Bias & $\sqrt{{\rm MSE}(\lambda)}$ & $\min\hat\lambda$ & $\max\hat\lambda$ & \% valid & $n$ \\
\hline
0.10 & 0.018 & 0.107 & 0.00 & 1.01 &  96 &  50 \\
 & 0.042 & 0.157 & 0.00 & 1.12 &  96 & 100 \\
 & 0.301 & 0.609 & 0.00 & 3.40 &  96 & 500 \\
 & 0.634 & 1.189 & 0.00 & 6.48 &  96 & 1000 \\
 0.25 & 0.009 & 0.132 & 0.00 & 0.91 &  99 &  50 \\
 & 0.006 & 0.161 & 0.00 & 0.95 &  99 & 100 \\
 & 0.123 & 0.427 & 0.00 & 2.27 &  99 & 500 \\
 & 0.320 & 0.776 & 0.00 & 4.00 &  99 & 1000 \\
 0.50 & 0.022 & 0.185 & 0.00 & 1.64 & 100 &  50 \\
 & 0.010 & 0.186 & 0.00 & 1.27 & 100 & 100 \\
 & 0.031 & 0.395 & 0.00 & 2.03 & 100 & 500 \\
 & 0.139 & 0.635 & 0.00 & 3.33 & 100 & 1000 \\
 0.75 & 0.031 & 0.239 & 0.08 & 1.99 & 100 &  50 \\
 & 0.014 & 0.215 & 0.00 & 1.90 & 100 & 100 \\
 & 0.002 & 0.394 & 0.00 & 2.25 & 100 & 500 \\
 & 0.049 & 0.612 & 0.00 & 3.28 & 100 & 1000 \\
1.00 & 0.045 & 0.304 & 0.31 & 2.68 & 100 &  50 \\
 & 0.021 & 0.256 & 0.26 & 2.24 & 100 & 100 \\
 & $-$0.000 & 0.388 & 0.00 & 2.65 & 100 & 500 \\
 & 0.012 & 0.608 & 0.00 & 3.69 & 100 & 1000 \\
1.50 & 0.061 & 0.419 & 0.57 & 4.22 & 100 &  50 \\
 & 0.028 & 0.330 & 0.59 & 3.22 & 100 & 100 \\
 & $-$0.004 & 0.380 & 0.03 & 3.00 & 100 & 500 \\
 & $-$0.014 & 0.580 & 0.00 & 3.97 & 100 & 1000 \\
 2.00 & 0.094 & 0.529 & 0.80 & 6.49 & 100 &  50 \\
 & 0.038 & 0.402 & 0.93 & 4.16 & 100 & 100 \\
 & 0.005 & 0.380 & 0.64 & 3.46 & 100 & 500 \\
 & $-$0.002 & 0.536 & 0.00 & 4.08 & 100 & 1000 \\
\hline
\end{tabular}
\end{center}
\caption{Empirical performance of the estimator $\dot\lambda_n$ of \eqref{eq:lamdot} given the estimate $\hat\sigma_n$ (see Table \ref{tab4}), for
different values of the parameter $\lambda$ and different sample
size. The time horizon $T$ is fixed to 500. Results over 10000
Monte Carlo paths of the geometric telegraph process. val = \% of valid cases, i.e.  simulated paths such that $\mu > \bar Y_n/\Delta$. See text for more details.} 
\label{tab3}
\end{table}

\begin{table}[ht]
\begin{center}
\begin{tabular}{c|cccc|cc}
\hline
 $\lambda$ & Bias & $\sqrt{{\rm MSE}(\sigma)}$ & $\min\hat\sigma_n$ & $\max\hat\sigma_n$ & \% valid & $n$ \\
&&&&& cases &\\
\hline
0.10 & $-$0.008 & 0.138 & 0.01 & 0.85 &  96 &  50, 100, 500, 1000 \\
% & $-$0.008 & 0.138 & 0.01 & 0.85 &  96 & 100 \\
% & $-$0.008 & 0.138 & 0.01 & 0.85 &  96 & 500 \\
% & $-$0.008 & 0.138 & 0.01 & 0.85 &  96 & 1000 \\
0.25 & $-$0.007 & 0.094 & 0.02 & 0.74 &  99 &  50, 100, 500, 1000 \\
% & $-$0.007 & 0.094 & 0.02 & 0.74 &  99 & 100 \\
% & $-$0.007 & 0.094 & 0.02 & 0.74 &  99 & 500 \\
% & $-$0.007 & 0.094 & 0.02 & 0.74 &  99 & 1000 \\
 0.50 & $-$0.003 & 0.065 & 0.09 & 0.69 & 100 &  50, 100, 500, 1000 \\
% & $-$0.003 & 0.065 & 0.09 & 0.69 & 100 & 100 \\
% & $-$0.003 & 0.065 & 0.09 & 0.69 & 100 & 500 \\
% & $-$0.003 & 0.065 & 0.09 & 0.69 & 100 & 1000 \\
0.75 & $-$0.002 & 0.053 & 0.22 & 0.66 & 100 &  50, 100, 500, 1000 \\
% & $-$0.002 & 0.053 & 0.22 & 0.66 & 100 & 100 \\
% & $-$0.002 & 0.053 & 0.22 & 0.66 & 100 & 500 \\
% & $-$0.002 & 0.053 & 0.22 & 0.66 & 100 & 1000 \\
1.00 & $-$0.002 & 0.045 & 0.30 & 0.66 & 100 &  50, 100, 500, 1000 \\
% & $-$0.002 & 0.045 & 0.30 & 0.66 & 100 & 100 \\
% & $-$0.002 & 0.045 & 0.30 & 0.66 & 100 & 500 \\
% & $-$0.002 & 0.045 & 0.30 & 0.66 & 100 & 1000 \\
1.50 & $-$0.002 & 0.037 & 0.33 & 0.64 & 100 &  50, 100, 500, 1000 \\
%& $-$0.002 & 0.037 & 0.33 & 0.64 & 100 & 100 \\
% & $-$0.002 & 0.037 & 0.33 & 0.64 & 100 & 500 \\
% & $-$0.002 & 0.037 & 0.33 & 0.64 & 100 & 1000 \\
2.00 & $-$0.001 & 0.031 & 0.36 & 0.61 & 100 &  50, 100, 500, 1000 \\
% & $-$0.001 & 0.031 & 0.36 & 0.61 & 100 & 100 \\
% & $-$0.001 & 0.031 & 0.36 & 0.61 & 100 & 500 \\
% & $-$0.001 & 0.031 & 0.36 & 0.61 & 100 & 1000 \\
\hline
\end{tabular}
\end{center}
\caption{Empirical performance of the estimator $\hat\sigma_n$ of \eqref{eq:sigma} for
different values of the parameter $\lambda$,  different sample
size and $\sigma=0.5$. The time horizon $T$ is fixed to 500. Results over 10000
Monte Carlo paths of the geometric telegraph process. `\% valid cases' = percentage of valid cases, i.e.  simulated paths such that the estimator of $\sigma$ exists, i.e. $\mu > \bar Y_n/\Delta$. See text for more details.} 
\label{tab4}
\end{table}

\section{Conclusions}
Despite the analytical problems in dealing with the telegraph process, in this paper we have shown that statistical inference on the process can be attempted. Moreover, numerical results seem to suggest that also for the geometric telegraph process this attempt might be successful which encourages the study of the analytical properties of this process in view of financial applications. In particular, a detailed description of the law of the increments of the process of the log-returns might be of interest because, numerical evidence, show that these have heavy tails but their law do not need a high number of parameters like other distributions proposed in the literature (see e.g. Eberlain and Keller, 1995). Moreover, the parameters have a direct interpretation as in the standard geometric Brownian motion.
This will be a topic for future research.

\end{document}